\documentclass[10pt]{article}
\pdfoutput=1
\usepackage{geometry, amsmath, amsthm, latexsym, amssymb, graphicx, mathtools, tikz, enumerate, dsfont}
\usepackage{tikz-cd}
\usepackage[utf8]{inputenc}
\usepackage[OT2,T1]{fontenc}
\usepackage[shortlabels]{enumitem}
\setenumerate{topsep=0pt}

\usepackage{tkz-graph}
\usepackage{graphicx}
\usepackage{csquotes}
\usepackage{lmodern}
\usepackage{arydshln}
\usepackage[shortlabels]{enumitem}
\usepackage{multirow}
\usepackage{mathdots}
\usepackage{hyperref}
\usepackage{authblk}
\usepackage{arydshln}
\usepackage{bigints}

\usepackage{mathrsfs}

\DeclareSymbolFont{cyrletters}{OT2}{wncyr}{m}{n}
\DeclareMathSymbol{\Sha}{\mathalpha}{cyrletters}{"58}

\theoremstyle{definition}
\newtheorem{defi}{Definition}[section]

\newtheorem{obs}[defi]{Remark}

\theoremstyle{plain}
\newtheorem{theorem}[defi]{Theorem}
\newtheorem{prop}[defi]{Proposition}
\newtheorem{lemma}[defi]{Lemma}
\newtheorem{cor}[defi]{Corollary}

\newtheorem*{idea*}{Idea}

\DeclareMathOperator{\Tr}{Tr}

\DeclareMathOperator{\SL}{SL}
\DeclareMathOperator{\GL}{GL}

\DeclareMathOperator{\So}{SO}

\DeclareMathOperator{\sgn}{sgn}
\DeclareMathOperator{\sgnrk}{sgnrk}
\DeclareMathOperator{\Disc}{Disc}

\DeclareMathOperator{\im}{Im}

\DeclareMathOperator{\vol}{vol}
\DeclareMathOperator{\covol}{covol}

\newcommand{\lp}{\left(}
\newcommand{\rp}{\right)}

\newcommand{\Z}{\mathbb{Z}}
\newcommand{\Q}{\mathbb{Q}}
\newcommand{\R}{\mathbb{R}}

\newcommand{\F}{\mathbb{F}}

\newcommand{\bH}{\mathbb{H}}
\newcommand{\pP}{\mathbb{P}}
\newcommand{\oO}{\mathcal{O}}

\newcommand{\fs}{\mathfrak{s}}

\newcommand{\cF}{\mathcal{F}}

\newcommand{\cK}{\mathcal{K}}

\newcommand{\eps}{\varepsilon}

\usepackage[
backend=biber,
style=alphabetic,
sorting=nyt
]{biblatex}
\title{Lower bounds on ranks of universal quadratic forms through geometry-of-numbers}
\author{Martí Oller}
\affil{\small Charles University, Faculty of Mathematics and Physics, Department of Algebra, Sokolovská 83, 186 75 Praha 8, Czech Republic}

\begin{document}
\maketitle
\begin{abstract}
We obtain lower bounds for the ranks of universal forms over totally
real cubic fields using techniques from geometry-of-numbers and the parametrisation
of cubic rings by binary cubics due to Delone--Faddeev. Additionally,
our methods give lower bounds for totally real primitive fields
of any degree with codifferent generated by a totally positive element.
\end{abstract}

\section{Introduction}

The study of quadratic forms and the integers they represent has a long
history, dating back essentially to the origins of mathematics itself.
In this article, we will focus on positive definite \emph{universal}
forms, that is, positive definite quadratic forms that represent all positive integers.
Over $\Z$, Lagrange proved that every positive integer can be represented
by the sum of four squares, or in other words that the form 
$x^2 + y^2 + z^2 + w^2$ is universal.

The analogous question over number fields is far from understood. 
Given a totally real number field $K$, we say
that a quadratic form $Q$ is \emph{totally positive definite} if all of its values
are totally positive (except for $Q(0,\dots,0) = 0$), and \emph{universal} if
it represents all totally positive elements of the ring of integers.
For instance, if $K = \Q(\sqrt{5})$ then Maass \cite{Maass} proved that
the sum of three squares is universal, and later Siegel \cite{Siegel}
showed that the only fields where a sum of squares is universal are
$\Q$ and $\Q(\sqrt{5})$, suggesting that universal forms of
low rank might be quite rare. In this direction, an influential conjecture by Kitaoka
asserts that only finitely many totally real number fields admit a universal
ternary quadratic form; this conjecture remains open.

One can consider these questions in the more general context of quadratic
lattices, which we now define explicitly.
Let $K$ be a totally real number field. We will consider positive definite
quadratic $\oO_K$-lattices $(\Lambda, Q)$, that is, finitely generated
$\oO_K$-modules $\Lambda$ equipped with a quadratic form $Q$ such that
for any non-zero $v \in \Lambda$ we have $Q(v) \in \oO_K^+$  (i.e. $Q(v)$
is totally positive). We will call such a lattice:
\begin{itemize}
\item \emph{universal} if for all $w \in \oO_K^+$ there exists $v \in \Lambda$
such that $Q(v) = w$.
\item \emph{classical} if all the values of the associated bilinear form
$B$ of $Q$ lie in $\oO_K$.
\end{itemize}
We will denote by $R(K)$ the smallest rank of a universal lattice over $K$, and
by $R_{cls}(K)$ the smallest rank of a universal classical lattice over
$K$.

For some specific cases, we have information about the behaviour of $R(K)$
and $R_{cls}(K)$. It is known that only finitely many quadratic fields
have $R(K) \leq 7$ \cite{KKP}, while there are infinitely many real
quadratic fields with $R_{cls}(K) = R(K) = 8$ \cite{Kim}. However, such
results appear to be rare, as Kala--Yatsyna--Żmija \cite{KYZ} proved that for
almost all real quadratic fields we have $R_{cls}(K) \gg_{\eps} D^{1/12-\eps}$
and $R(K) \gg_{\eps} D^{1/24-\eps}$, where $D$ denotes the discriminant of the field and $\eps > 0$.
Less is known in higher degree: some analogous results are known for
multiquadratic fields \cite{sails, Man} or for specific families of
cubic fields \cite{KT}.

In this article, we obtain lower bounds for $R(K)$ and $R_{cls}(K)$
for most totally real cubic fields (conjecturally $98\%$). Our methods are general enough
that we can obtain comparable lower bounds for other families of primitive
fields (i.e. fields $K/\Q$ with no intermediate fields) in arbitrarily
high degree.

Given a totally real number field $K$ of degree $d$, denote its 
distinct real embeddings by $\sigma_1,\dots,\sigma_d \colon K \xhookrightarrow{} \R$. 
Let $\alpha \in K$: its \emph{signature} is the $d$-tuple
\[
\fs(\alpha) = (\sgn(\sigma_1(\alpha)),\dots,\sgn(\sigma_d(\alpha))) \in \{\pm\}^d.
\]
If $\fs(\alpha) = (+,\dots,+)$, we call $\alpha$ \emph{totally positive}.
We denote by $K^+$ the set of totally positive elements of $K$,
and we further denote $\oO_K^+ = \oO_K \cap K^+$.
Consider the $\F_2$-vector space $\fs(\oO_K^{\times})$ of signatures attained
by the units $\oO_K^{\times}$: we call its dimension the \emph{signature rank},
denoted $\sgnrk(K)$. We note that $1 \leq \sgnrk(K) \leq d$, as $-1$ is always
a unit.

\begin{theorem}
\label{theo:sgnrk3}
Let $K$ be a totally real cubic field with $\sgnrk(K) = 3$ and discriminant
$D$. Then
\[
R_{cls}(K) \geq \frac{\sqrt{7}}{12} D^{1/4} - \frac{5}{12} \quad
\text{and} \quad 
R(K) \geq \begin{cases}
\frac{\sqrt[4]{7}}{6\sqrt{2}} D^{1/8} - \frac{\sqrt{5}}{6\sqrt{2}} & \text{if } D < D_0, \\
\frac{\sqrt[4]{7}}{6} D^{1/8} - \frac{\sqrt{5}}{6} & \text{if } D \geq D_0,
\end{cases}
\]
where $D_0 = 17\,405\,939\,540$.
\end{theorem}

\begin{theorem}
\label{theo:sgnrk2}
Let $\psi \colon \R_{\geq 0} \to \R_{\geq 0}$ be a function with $\lim_{D \to \infty} \psi(D) = +\infty$.
Then, when totally real cubic fields $K$ with $\sgnrk(K) = 2$ are ordered
by discriminant, $100 \%$ of them satisfy
\[
R_{cls}(K) \geq \frac{D^{1/4}}{\psi(D)} \quad \text{and} \quad R(K) \geq \frac{D^{1/8}}{\psi(D)}.
\]
\end{theorem}
We note that by \cite[Corollary 10]{BVarma}, when ordered by discriminant at least $75\%$ of totally real cubic fields
have signature rank $\geq 2$. Dummit--Voight conjectured in \cite{DV} that the approximate
densities for each signature rank when ordered by discriminant are
$36.3 \%$ for $\sgnrk(K) = 3$, $61.8 \%$ for $\sgnrk(K) = 2$ and $1.9 \%$ for $\sgnrk(K) = 1$.
So, if the conjectural densities hold, then Theorems \ref{theo:sgnrk3}
and \ref{theo:sgnrk2} cover approximately $98 \%$ of all totally real cubic fields.

We also obtain results for certain fields of arbitrary degree. Recall that a 
field extension $K/F$ is \emph{primitive} if there are no fields $L$ such that
$F \subsetneq L \subsetneq K$. If $K$ is a number field, we say it is
\emph{primitive} if the extension $K/\Q$ is primitive.
\begin{theorem}
\label{theo:nu1}
Let $K/\Q$ be a totally real primitive number field of degree $d$, and assume that
$\oO_K^{\vee} = (\delta)$ for some $\delta \in K^+$. Then
\[
R_{cls}(K) \geq \frac{D^{1/2}}{d!} + O_d(D^{\frac{1}{2}-\frac{1}{d(d-1)}})
\]
and
\[
R(K) \geq \frac{D^{1/4}}{d\sqrt{(d-1)!}} + O_d(D^{\frac{1}{4}-\frac{1}{2d(d-1)}}).
\]
The implicit constant in the error term depends on $d$, but not on $K$.
\end{theorem}

In fact, we get a slightly more general result in Corollary \ref{cor:counting}.
Theorem \ref{theo:nu1} generalises a result by Yatsyna \cite[Theorem 4]{Yatsyna} by
weakening the hypotheses and giving an explicit order of growth, and it is
a variant of a result by Kala--Melistas \cite[Theorem 1.2]{KM}, which
does not contain the primitivity assumption but also does not give an
explicit order of growth.

The main idea behind our proofs is Proposition \ref{prop:trace1}, a 
result by Yatsyna \cite[Theorem 25]{Yatsyna} that shows that if there
exists $\delta \in \oO_K^{\vee,+}$ such that there are many elements
$\alpha_i \in \oO_K^+$ with $\Tr(\delta\alpha_i) = 1$, then $R(K)$ and
$R_{cls}(K)$ are forced to be large. With this in mind, in Section \ref{section:gon}
we use geometry-of-numbers methods to obtain an estimate on how many such
$\alpha_i$ might exist for a given $\delta$ for primitive fields. Then,
for cubic fields we use the Delone--Faddeev parametrisation (recalled in
Section \ref{section:DF}) to find suitable candidates $\delta$ of
low enough norm: this is carried out in Section \ref{section:mainPfs}.

\subsection{Acknowledgements}

The author would like to thank Vítězslav Kala, Siu Hang Man and Pavlo
Yatsyna for helpful discussions around the paper.
The article was completed while the author was supported by the Czech 
Science Foundation GAČR, grant 26-20514S.

\subsection{AI disclosure}

This article was prepared with the help of Claude, model Opus 5.5. The
main underlying approach is due to the human author, and AI was used
to generate the proof of Lemma \ref{lemma:width} and the mathematical
content of Sections \ref{subs:sgnrk3} and \ref{subs:sgnrk2}. For all these
results, the AI output was modified (sometimes significantly) to be 
presented in a more readable way. Every word in this article is 
human-written, and AI was used for incidental proof-writing 
(mathematical and non-mathematical).

\section{Counting lattice points in codifferent-trace hyperplanes}
\label{section:gon}

Let $K$ be a totally real number field. Its codifferent is defined as
\[
\oO_K^{\vee} := \{\alpha \in K \mid \Tr(\alpha\beta) \in \Z, \, \forall \beta \in \oO_K\}.
\]
Given $\delta \in \oO_K^{\vee,+}$, let us define
\[
S(\delta) = \{\alpha \in \oO_K^+ \mid \Tr(\delta \alpha) = 1\}.
\]

\begin{prop}
\label{prop:Sdelta}
Let $d \geq 2$, and let $K$ be a totally real number field of degree $d$
of discriminant $D$. Let $\delta \in \oO_K^{\vee,+}$ 
be a primitive element, that is, such that if $\delta/m \in \oO_K^{\vee,+}$
for some integer $m \geq 1$, then $m = 1$. Then
\[
\# S(\delta) = \frac{\sqrt{D}}{(d-1)! \nu} + O_d \lp \frac{D^{\frac{1}{2}-\frac{1}{d(d-1)}}}{\nu^{1-2/d}} \rp,
\]
where $\nu = N(\delta)D$. The implied constant depends on $d$, but not
on $K$ or $\delta$.
\end{prop}
Note that $\nu \in \Z_{\geq 1}$.
The main technical tool in this proof is Davenport's lemma, proven in
\cite{Davenport}, for which we state a particular version below.

\begin{lemma}[Davenport's lemma]
\label{lemma:Davenport}
Let $S \subset \R^n$ be a convex, bounded set. For $0 \leq m \leq n-1$,
denote by $S_m$ the sum of the volumes of the projections of $S$ on the $m$-dimensional
coordinate hyperplanes obtained by equating any $n-m$ coordinates to zero.
Then
\[
|\#(S \cap \Z^n) - \vol(S)| \leq \sum_{m=0}^{n-1} S_m.
\]
\end{lemma}

\begin{obs}
Lemma \ref{lemma:Davenport} corresponds to taking $h = 1$ in \cite[Theorem]{Davenport},
which we can do since we are additionally assuming that $S$ is convex.
\end{obs}

Now, consider a totally real number field $K$ of degree $d$, together
with the Minkowski embedding $\sigma \colon K \xhookrightarrow{} \R^d$.
Then, $\sigma(\oO_K)$ (which we will write 
simply as $\oO_K$ by abuse of notation) is a lattice in $\R^d$ of rank 
$d$. Let $\delta \in \oO_K^{\vee,+}$ be primitive in the sense of 
Proposition \ref{prop:Sdelta}, and consider the hyperplane $H_{\delta} \subset \R^d$
defined by $H_{\delta} := \{x \in \R^d \mid \Tr(\delta x) = 1\}$. We note that primitivity
implies that $H_{\delta} \cap \oO_K \neq \emptyset$, and thus that $H_{\delta} \cap \oO_K$
is a full-rank lattice inside $H_{\delta}$. To determine the main term in Davenport's
lemma, first we need to determine the volume of $T_{\delta} := H_{\delta} \cap (\R_{>0})^d$, normalised by the covolume
of the lattice $\Gamma := H_{\delta} \cap \oO_K$.

Note that $T_{\delta}$ intersects the coordinate axes
at $e_i/\delta_i$, where $e_i$ is the $i$-th unit vector and $\delta = (\delta_1,\dots,\delta_d)$
are the coordinates of $\delta$ under the Minkowski embedding. The volume
of $T_{\delta}$ is obtained through the following:
\begin{lemma}
\label{lemma:nSimplex}
Let $v_0,\dots,v_n$ be vertices of an $n$-simplex inside $\R^n$, and let
$A$ be the matrix with columns $v_i - v_0$ for all $i = 1,\dots,n$. Then
the $n$-dimensional volume of the $n$-simplex is equal to $\frac{1}{n!}|\det(A)|$.
\end{lemma}
Through Lemma \ref{lemma:nSimplex}, the volume of $T_{\delta}$ can be computed to be
\[
\vol(T_{\delta}) = \frac{\sqrt{\delta_1^2 + \dots + \delta_d^2}}{(d-1)! \delta_1\dots\delta_d}.
\]
We note that $\delta_1\dots \delta_d = N(\delta)$. It remains to determine
the covolume of the lattice $\Gamma = H_{\delta} \cap \oO_K$ inside $H_{\delta}$.
We note that the different consecutive hyperplanes $\{x \in \R^d \mid \Tr(\delta x) = n\}$ 
for $n \in \Z$ are at a distance $1/\sqrt{\delta_1^2 + \dots + \delta_d^2}$
to each other, so we get that
\[
\covol(\Gamma) \frac{1}{\sqrt{\delta_1^2 + \dots + \delta_d^2}} = \covol(\oO_K) = \sqrt{D},
\]
so $\covol(\Gamma) = \sqrt{D(\delta_1^2 + \dots + \delta_d^2)}$. In conclusion,
the normalised volume of $T_{\delta}$ is $\frac{1}{(d-1)! N(\delta) \sqrt{D}}$.
Writing $\nu = N(\delta)D$, we can rewrite this volume as 
\begin{equation}
\label{eq:volT}
\vol_{\Gamma}(T_{\delta}) = \frac{\sqrt{D}}{(d-1)! \nu}.
\end{equation}
Now, it remains to control how big the projections of $T_{\delta}$ onto 
coordinate hyperplanes can be taken to be. A key tool for us is the
realisation that we can choose a convenient basis for the lattice $\Gamma$
in such a way that the region $T_{\delta}$ is not very skewed.
\begin{defi}
Let $L \subset \R^n$ be a lattice with inner product $\langle \cdot, \cdot \rangle$, 
and let $S \subset \R^n$ be a closed,
bounded, convex set. The \emph{width} $w_L(S)$ of $S$ with respect to $L$
is defined as
\[
w_L(S) := \min_{\alpha \in L^{\vee}\setminus\{0\}} \lp \max_{x \in S}\langle \alpha,x\rangle - \min_{y \in S}\langle \alpha,y\rangle\rp,
\]
where $L^{\vee}$ denotes the dual lattice of $L$.
\end{defi}

When $K$ is primitive, we can obtain a lower bound on $w_{\Gamma}(T_{\delta})$:
\begin{lemma}
\label{lemma:width}
Suppose that $K$ is a totally real primitive number field of degree $d$.
Then there exists a constant $c_d > 0$, depending only on $d$, such that
\[
w_{\Gamma}(T_{\delta}) \geq c_d D^{1/(d(d-1))} \nu^{-2/d}.
\]
\end{lemma}
\begin{proof}
We note that every element of $\Gamma^{\vee}$ is the restriction to
$H_{\delta}$ of $\Tr(\alpha \cdot)$ of some $\alpha \in \oO_K^{\vee}$,
and this restriction is non-zero if $\alpha \notin \Q\delta$ (i.e. those
directions perpendicular to $\delta$). Hence, to obtain a lower bound we can
consider all the lattice width directions corresponding to elements
$\alpha \in \oO_K^{\vee}$ in $\R^d$ instead of just $\alpha \in \Gamma^{\vee}$
in $H_{\delta}$.

Let $\alpha \in \oO_K^{\vee}$. Then the maximum and the minimum $\Tr(\alpha x)$
for $x \in T_{\delta}$ are always attained by one of the vertices of $T_{\delta}$, which
we recall are of the form $e_i/\delta_i$. Then the problem reduces to finding
\[
\min_{\rho}\max_{1 \leq i < j \leq d}|\rho_i - \rho_j|
\]
among all $\rho \in \delta^{-1} \oO_K^{\vee}$ such that $\rho \notin \Q$. 
We note that for any such $\rho$ we have $\Q(\rho) = K$ since $K$
is primitive. Hence $M := \Z \oplus \Z\rho \oplus \dots \oplus \Z\rho^{d-1}$
is a lattice of full rank in $\R^d$, and $M \subset (\delta^{-1} \oO_K^{\vee})^{d-1}$
since $\rho^{d-1} \in (\delta^{-1} \oO_K^{\vee})^{d-1}$. Thus, the covolume of
$M$ in the Minkowski embedding is
\[
\covol(M) \geq N(\delta)^{-(d-1)} D^{-(2d-3)/2} = D^{1/2}\nu^{-(d-1)}.
\]
On the other hand $\covol(M) = \prod_{i < j} |\rho_i-\rho_j|$.
If $\rho_1 \leq \rho_2 \leq \dots \leq \rho_d$ then there exists some constant
$c_d' \leq 1$ such that
\[
\covol(M) \leq c_d' |\rho_d - \rho_1|^{d(d-1)/2}.
\]
For instance, we can take $c_3' = 1/4$, and $c_d = 1$ always works.
Putting it all together, we conclude that
\[
w_{\Gamma}(T_{\delta}) \geq c_d D^{1/(d(d-1))} \nu^{-2/d},
\]
for $c_d = (c_d')^{-2/(d(d-1))}$, as wanted.
\end{proof}
\begin{obs}
\label{obs:c3}
For any $d$, the lemma holds with $c_d = 1$. For $d = 3$ we can take
$c_3 = \sqrt[3]{4}$.
\end{obs}

Proposition \ref{prop:Sdelta} now follows from combining Lemma \ref{lemma:Davenport},
Equation \eqref{eq:volT} and Lemma \ref{lemma:width}.

\subsection{Refinements for cubic sails}

If $[K \colon \Q] = 3$, then $K$ is always primitive and Proposition \ref{prop:Sdelta}
always applies. We can in fact obtain explicit constants on the implicit
error terms in the proposition; however, to obtain our main results we
will need some more precise results.

\begin{lemma}
\label{lemma:Bohnert}
Let $S \subset \R^2$ be a compact, convex subset with non-empty interior,
lattice width $w > 1$ and area $A$. Then
\[
\# (S \cap \Z^2) \geq \begin{cases}
\frac{2(w-1)}{w^2}A - 1 & \text{if }1 < w \leq 2, \\
\frac{1}{2}A - \frac{5}{4} & \text{if } w > 2.
\end{cases}
\]
\end{lemma}
\begin{proof}
This is a combination of \cite[Theorems 1.1, 1.2]{Bohnert} in the case when
$\# (S \cap \Z^2) \geq 1$ and \cite[Theorem 2.2]{AW} when $\# (S \cap \Z^2) = 0$.
\end{proof}

With this result, the leading constant in the main term becomes worse
but the error term becomes essentially constant, as long as we can guarantee
that the lattice width is strictly greater than $1$. Given $\delta \in \oO_K^{\vee,+}$
and the triangle $T$ defined by the intersection of the hyperplane $H \colon \Tr(\delta \cdot) = 1$
with the totally positive octant, we know from the last section that the
lattice width of $T$ with respect to $\Gamma := H \cap \oO_K$ is at least
$\sqrt[3]{4} D^{1/6} \nu^{-2/3}$. This is guaranteed to be larger than
$1$ whenever $\nu < 2D^{1/4}$.

\section{The Delone--Faddeev correspondence}
\label{section:DF}

A cubic ring is a (commutative, unitary) ring that is also a free $\Z$-module of rank $3$. By
a theorem of Delone--Faddeev \cite{DF}, later strengthened by Gan--Gross--Savin
\cite{GGS}, these cubic rings are parametrised by equivalence classes
of binary cubic forms. Here, a binary cubic form $f(x,y) = ax^3 + bx^2y + cxy^2 + dy^3$
is acted upon by $\gamma \in \GL_2(\Z)$ via the twisted action
\[
(\gamma \cdot f)(x,y) = \frac{1}{\det \gamma}f((x,y)\cdot \gamma).
\]
\begin{theorem}
\label{theorem:D-F}
There is a canonical, discriminant-preserving bijection between isomorphism classes of cubic rings
and $\GL_2(\Z)$-equivalence classes of integral binary cubic forms.
\end{theorem}
This bijection can be made very explicit, as can be seen in \cite[\S 2]{BST},
for instance. In \cite{Wood}, Wood establishes a more general equivalence
than Theorem \ref{theorem:D-F} for binary $n$-ic forms, in particular
finding explicit $\Z$-bases for the ring of integers and the codifferent.
We now write these bases explicitly.

Let $K$ be a cubic number field (totally real or not), and consider its
ring of integers $\oO_K$. By Theorem \ref{theorem:D-F}, it corresponds
to a binary cubic $f(x,y) = ax^3 + bx^2y + cxy^2 + dy^3$. If
$K_f = \Q[x]/(f(x,1))$ and $\xi$ is the image of $x$ in $K_f$, then
we have that
\[
\oO_K \cong \langle 1, a\xi, a\xi^2 + b\xi\rangle.
\]
With this, a basis for the codifferent can be explicitly determined
by \cite[Theorem 2.2, Corollary 3.7]{Wood} to be
\[
\oO_K^{\vee} \cong \frac{1}{f'(\xi,1)} \langle 1, \xi, a\xi^2\rangle.
\]
Note that if the binary cubic has $x^3$-coefficient equal to
$\pm 1$, then $\oO_K$ is monogenic and we recover the usual expression
for the codifferent in monogenic fields $\oO_K^{\vee} = 1/f'(\xi) \oO_K$,
if $\oO_K = \Z[\xi]$.

\section{Proof of the main theorems}
\label{section:mainPfs}

The main technical tool for proving lower bounds on $R_{cls}(K)$ and $R(K)$
is the following result by Yatsyna \cite[Theorem 25]{Yatsyna}, cf.
\cite[Propositions 7.2 and 7.4]{KT}.

\begin{prop}
\label{prop:trace1}
Let $K$ be a totally real number field of degree $d$ over $\Q$, and $\oO$
an order in $K$. Suppose that there exist
\begin{itemize}
\item an element $\delta \in \oO^{\vee,+}$,
\item $n$ elements $\alpha_1,\dots,\alpha_n \in \oO^+$ such that $\Tr(\delta \alpha_i) = 1$ for all $i$.
\end{itemize}
Then
\[
R_{cls}(\oO) \geq \frac{n}{d} \quad \text{ and } \quad R(\oO) \geq \begin{cases}
\frac{\sqrt{n}}{d} & \text{if } n \geq 240, \\
\frac{\sqrt{n}}{d\sqrt{2}} & \text{if } n < 240.
\end{cases}
\]
\end{prop}

Thus, in view of Propositions \ref{prop:Sdelta} and \ref{prop:trace1}, 
it suffices to find $\delta \in \oO_K^{\vee,+}$ with
a low enough $\nu$ relative to the discriminant to obtain lower bounds
for $R(K)$ and $R_{cls}(K)$. Precisely:
\begin{cor}
\label{cor:counting}
Let $K$ be a totally real primitive field of degree $d$ and discriminant
$D$, and let $\nu = \min_{\delta \in \oO_K^{\vee,+}} N(\delta) D$. Then
\[
R_{cls}(K) \geq \frac{\sqrt{D}}{d! \nu} + O_d \lp \frac{D^{\frac{1}{2}-\frac{1}{d(d-1)}}}{\nu^{1-2/d}} \rp
\]
and
\[
R(K) \geq \frac{D^{1/4}}{d \sqrt{(d-1)! \nu}} + O_d \lp \frac{D^{\frac{1}{4}-\frac{1}{2d(d-1)}}}{\nu^{\frac{1}{2}-\frac{1}{d}}} \rp.
\]
The implicit constants depend on $d$ but not on $K$.
\end{cor}

In particular, Theorem \ref{theo:nu1} follows directly from Corollary
\ref{cor:counting}, since in that context we always have $\nu = 1$.
More generally, we get non-trivial bounds for families of totally real
primitive degree $d$ fields with $\nu = o(D^{1/(2(d-1))})$.

For the rest of the section, we will focus on the cubic case and how the
different signature rank conditions might guarantee a choice of $\delta$
with an appropriately small norm.

\subsection{The case $\sgnrk(K) = 3$}
\label{subs:sgnrk3}

We will now prove Theorem \ref{theo:sgnrk3} for cubic fields $K$ with
signature rank $3$, that is, such that there are units of all possible
signatures, which will give many possible options to obtain an element
$\delta \in \oO_K^{\vee,+}$ of low norm.

\begin{lemma}
Let $K$ be a totally real cubic number field with $\sgnrk(K) = 3$ and
discriminant $D$. Then there always exists $\delta \in \oO_K^{\vee,+}$ 
such that $N(\delta) D \leq \lp \frac{D}{49} \rp^{1/4}$.
\end{lemma}
\begin{proof}
Consider the binary cubic form $f(x,y) = ax^3 + bx^2y + cxy^2 + dy^3$
associated to $\oO_K$ by the Delone--Faddeev correspondence. In particular,
recall that inside the field $\Q[x]/(f(x,1))$, if $\xi$ denotes a root of $f(x,1)$
we can identify
\[
\oO_K = \langle 1, a\xi, a\xi^2 + b\xi \rangle, \quad \oO_K^{\vee} = \frac{1}{f'(\xi,1)} \langle 1, \xi, a\xi^2 \rangle
\]
Consider $P = (p \colon q) \in \pP^1(\Q)$, where $p,q$ are integers, not
both simultaneously zero. The element $\delta_P = (q\xi - p)/f'(\xi,1)$ 
is then an element of the codifferent, and it can be shown that the norm
is $|N(\delta_P)| = |f(p,q)|/D$. By a theorem of Mordell \cite{Mordell},
we can choose $P = (p,q) \in \Z^2$, $(p,q) \neq (0,0)$, such that
\[
|f(p,q)| \leq \lp \frac{D}{49} \rp^{1/4}.
\]
We can multiply $\delta_P$ by a unit of appropriate signature, so 
we may assume that $\delta_P$ is totally positive, which concludes the proof.
\end{proof}

By Lemma \ref{lemma:width} and Remark \ref{obs:c3}, the lattice width of 
the triangle $T_{\delta}$ is larger than $\sqrt[3]{28} \approx 3.037$. 
Thus, we can use Lemma \ref{lemma:Bohnert} to conclude that
\[
\# S(\delta) \geq \frac{\sqrt{7}}{4} D^{1/4} - \frac{5}{4},
\]
and by Proposition \ref{prop:trace1} we get
\[
R_{cls}(K) \geq \frac{\sqrt{7}}{12} D^{1/4} - \frac{5}{12} \quad
\text{and} \quad 
R(K) \geq \frac{\sqrt[4]{7}}{6} D^{1/8} - \frac{\sqrt{5}}{6},
\]
where the inequality for $R(K)$ holds as long as $\# S(\delta) \geq 240$,
and there is an extra factor of $1/\sqrt{2}$ otherwise. This concludes the proof
of Theorem \ref{theo:sgnrk3}.

\subsection{The case $\sgnrk(K) = 2$}
\label{subs:sgnrk2}

When $\sgnrk(K) = 2$, the strategy of the proof is largely the same as when $\sgnrk(K) = 3$
but now we need to adapt for the fact that we can't choose $\delta_P$
so arbitrarily.

\begin{lemma}
Let $K$ be a totally real cubic field. Then there exist $p,q \in \Z$, not
both zero, such that $\delta_P = (q\xi - p)/f'(\xi,1)$ has the same
signature as some unit $u \in \oO_K^{\times}$.
\end{lemma}
\begin{proof}
Under a suitable ordering of the three real embeddings $\sigma_1,\sigma_2,\sigma_3$
of $K$, the signature of $\xi - p/q$ can be $(-,-,-), (+,-,-), (+,+,-)$ or $(+,+,+)$. We also note that
the signature of $f'(\xi,1)$ is, up to multiplication by $-1$, equal
to $(+,-,+)$. Thus, the set of elements of the form $(q \xi - p)/f'(\xi,1)$
attains the signatures $\pm(+,-,+), \pm (+,+,-)$ and $\pm (-,+,+)$.
Because $\sgnrk(K) = 2$, there has to exist a unit of $K$ with one
of those signatures.
\end{proof}

Let $t_1,t_2,t_3$ be the three roots of $f(x,1)$, where we recall that
$f$ is the binary cubic form from the Delone--Faddeev correspondence.
In practice, what we will do is to split $\pP^1(\R) - \{t_1,t_2,t_3\}$
into three different arcs $X_1,X_2,X_3$, where the signature of $\delta_P$
remains constant (up to multiplication by $-1$) for all $P \in X_i$, 
$i = 1,2,3$, and changes in the different arcs.
Let $P = (p \colon q) \in \pP^1(\Q)$, and consider $P' = (p' \colon q') \in \pP^1(\Q)$,
where $p,q,p',q'$ are integers with $|pq' - p'q| = 1$. Then, if
\[
M_{P,P'}(t) = \frac{q't-p'}{qt-p},
\]
the quantity $\max_{i=1,2,3}M_{P,P'}(t_i) - \min_{i=1,2,3}M_{P,P'}(t_i)$
is independent of $P'$, and we will denote it by $s_f(P)$.

\begin{lemma}
\label{lemma:sfP}
Let $f(x,y)$ be a binary cubic form with distinct real roots $t_1,t_2,t_3$. 
Suppose that $\pP^1(\R) - \{t_1,t_2,t_3\} = X_1 \sqcup X_2 \sqcup X_3$,
where the $X_i$ are connected.
Then, there exists $P = (p \colon q) \in X_i \cap \pP^1(\Q)$ such that
$s_f(P)  > 1$.
\end{lemma}
\begin{proof}
Suppose that one of the endpoints of $X_i$ is $t_i$, and consider
the convergents $x_n = (p_n,q_n)$ of the continued fraction of $t_i$.
For high enough $n$, choose $P = (p_n \colon q_n)$ and $P' = (p_{n-1} \colon q_{n-1})$,
where the parity of $n$ is chosen so that $P \in X_i$. It is an elementary
property of continued fractions that $p_nq_{n-1} - p_{n-1}q_n = (-1)^{n-1}$.

We note that $|M_{P,P'}(t_i)| > 1$, since the sequence $|q_k t_i - p_k|$
is strictly decreasing. We also note that $M_{P,P'}(t_j)$ has opposite
sign to $M_{P,P'}(t_i)$ if $i \neq j$, so in that case
\[
|M_{P,P'}(t_i) - M_{P,P'}(t_j)| > |M_{P,P'}(t_i)| > 1,
\]
as wanted.
\end{proof}

Having to consider the arcs separately means we can't apply Mordell's
theorem on the minimum attained by a binary cubic form, and so the constants
that we would obtain would not be good enough for a bound for \emph{every}
field $K$ with $\sgnrk(K) = 2$. Instead, we will obtain bounds that will
work everywhere except for a set of fields of density $\eps$, where we
will be able to make $\eps > 0$ arbitrarily small.

\begin{lemma}
\label{lemma:densityEps}
If $K$ is a totally real cubic field, denote by $f_K$ its associated
binary cubic form.
Let $\eps \in (0,1)$. Then there exists a compact set $S_{\eps}$ inside
the set of real binary cubic forms of discriminant $1$ (with the usual
Euclidean topology induced by $\R^{4}$) such that
\[
\limsup_{X \to \infty} \frac{\# \{K \mid \Disc(K) \leq X, \, (\Disc(K))^{-1/4}f_K \notin S_{\eps} \}}{\# \{K \mid \Disc(K) \leq X \}} \leq \eps.
\]
\end{lemma}
\begin{proof}
This result largely follows from the equidistribution of shapes of lattices
of real cubic fields, proven in \cite{BHshapes}. More precisely, following \cite{BHshapes},
we can associate to $\oO_K$ a \emph{shape}, that is, an element of $\SL_2(\Z) \backslash \bH$
where $\bH$ denotes the upper complex half-plane which is acted upon by
the usual Möbius action. Working out the map explicitly, the shape
associated to the ring of integers of a totally real cubic field $K$ is
given by a root $\tau$ of the Hessian of $f_K$ (which always has a pair of
complex roots if $f_K$ has three real roots).

Let $V_1$ denote the set of real binary cubic forms with discriminant
$1$, and let $\pi \colon V_1 \to \bH$ denote the map
sending $f$ to the root of its Hessian with positive imaginary part. 
The group $\SL_2(\R)$ acts transitively on $V_1$ via change of variables
and on $\bH$ via Möbius transformations, and $\pi$ is equivariant with
respect to this action. We also have that $\pi(v) = \pi(g \cdot v)$
if and only if $g$ is in the stabiliser of $v$, which is isomorphic to
a conjugate of $\So_2(\R)$, which is a compact subgroup. It follows that
if $C \subset \bH$ is compact, then $\pi^{-1}(C)$ is compact.

Take $\cF= \{\tau \in \bH \mid |\text{Re}(\tau)| \leq \frac{1}{2}, \, |\tau| \geq 1 \}$,
which is a fundamental domain for the action of $\SL_2(\Z)$ on $\bH$.
Consider the subset $C_Y := \{\tau \in \cF \mid \im \tau \leq Y\}$ for some sufficiently
large $Y$: under the hyperbolic measure $\mu$, we have that $\mu(\cF) < \infty$
and
\[
1-\frac{\mu(C_Y)}{\mu(\cF)} \leq \frac{3}{\pi Y},
\]
so we can take $Y$ such that $\frac{\mu(C_Y)}{\mu(\cF)} \geq 1-\eps$.

Take $S_{\eps} = \pi^{-1}(C_Y)$, a compact subset of $V_1$. By \cite[Theorem 1]{BHshapes},
the set of totally real cubic fields $K$ such that $D^{-1/4}f_K$ has an
$\SL_2(\Z)$-translate lying in $S_{\eps}$ has density $\frac{\mu(C_Y)}{\mu(\cF)}$,
which is greater than or equal to $1-\eps$, as wanted.
\end{proof}

\begin{theorem}
\label{theo:sgnrk2eps}
Let $\cK_{3,2}$ denote the set of totally real cubic fields with signature
rank $2$. Let $\eps \in (0,1)$. Then there exist constants $c_{\eps} > 0$ and $c_{\eps}'$
and a subset $\cK_{\eps} \subset \cK_{3,2}$ of density at least $1-\eps$ such that
for every $K \in \cK_{\eps}$ there exists $\delta \in \oO_K^{\vee,+}$
with
\[
\# S(\delta) \geq c_{\eps} D^{1/4} - c_{\eps}'.
\]
\end{theorem}
\begin{proof}
We first note that a positive density of totally real cubic fields
have signature rank equal to $2$ by \cite[Corollaries 9 and 10]{BVarma}. Thus, by Lemma \ref{lemma:densityEps}
we can choose a compact subset $S = S_{\eps'}$ (for some $\eps'$) inside $V_1$, the set of real binary
cubics of discriminant 1, such that for a density of at least $1-\eps$
of $K \in \cK_{3,2}$ with associated binary cubic form $f_K$, we have
that $D^{-1/4}f_K \in S$.

Take $f \in S$, and denote by $t_1,t_2,t_3$ the three distinct real roots
of $f(t,1)$. Then $\pP^1(\R) \setminus \{t_1,t_2,t_3\}$ can be split into
three arcs $X_1,X_2,X_3$, and for each arc $X_i$ there always exists $P_i \in X_i$
such that $s_f(P_i) > 1$ by Lemma \ref{lemma:sfP}. By continuity, there exists an open neighbourhood $U$
of $f$ in $S$ such that $P_i$ is not a root of $g$ and $s_{g}(P_i) > 1$
for all $g \in U$. By compactness, we can cover $S$ with finitely many
such neighbourhoods, meaning that there exist constants $C_1,C_2,C_3 > 0$ such that
for every $f \in S$ and every arc $X$ of $f$, there exists $P \in X$
such that
\[
1 + C_1 \leq|s_f(P)| \leq C_2, \quad |f(P)| \leq C_3.
\]
Let $K \in \cK_{3,2}$ be such that $D^{-1/4}f_K \in S$, choose an arc $X$ of $f_K$
such that there exists a unit of the corresponding signature and take $\delta_P$ with
$P$ as above. Then, the area of the triangle defined by $\delta_P \in \oO_K^{\vee,+}$ is
$A \geq D^{1/4}/(2C_3)$, and the lattice width greater than or equal
to $\min\{s_f(P), 2A/s_f(P)\}$, which is $>1$ for all but finitely many
fields $K$. The result then follows from Lemma \ref{lemma:Bohnert}.
\end{proof}

Finally, Theorem \ref{theo:sgnrk2} follows from Theorem \ref{theo:sgnrk2eps}.

\printbibliography

\end{document}